\documentclass[12pt]{article}
\usepackage[utf8]{inputenc}
\usepackage{fullpage}
\usepackage[sc]{mathpazo}
\usepackage[T1]{fontenc}
\usepackage[margin=10pt,font=small]{caption}

\usepackage{amsmath}
\usepackage{amssymb}
\usepackage{amsthm}
\usepackage{graphicx}
\usepackage{hyperref}
\usepackage{color}

\newtheorem{thm}{Theorem}[section]
\newtheorem{prop}[thm]{Proposition}
\newtheorem{lem}[thm]{Lemma}
\newtheorem{coro}[thm]{Corollary}
\theoremstyle{remark}
\newtheorem{rmk}{Remark}

\newcommand{\tdef}[1]{\textcolor{blue}{\emph{#1}}}

\newcommand{\type}{\mathrm{Type}}

\newcommand{\sdd}{\mathbf{c}_{00}}
\newcommand{\sdu}{\mathbf{c}_{01}}
\newcommand{\suu}{\mathbf{c}_{11}}
\newcommand{\rcont}{\mathbf{rcont}}

\newcommand{\black}{\mathbf{black}}
\newcommand{\white}{\mathbf{white}}
\newcommand{\face}{\mathbf{face}}
\newcommand{\outdeg}{\mathbf{outdeg}}

\newcommand{\lnode}{\mathbf{lnode}}
\newcommand{\znode}{\mathbf{znode}}
\newcommand{\pnode}{\mathbf{pnode}}
\newcommand{\rlabel}{\mathbf{rlabel}}

\newcommand{\maptotree}{\mathbf{T}_\mathcal{M}}
\newcommand{\treetomap}{\mathbf{M}_\mathcal{T}}
\newcommand{\treetoint}{\mathbf{I}_\mathcal{T}}
\newcommand{\inttotree}{\mathbf{T}_\mathcal{I}}
\newcommand{\maptoint}{\mathbf{I}_\mathcal{M}}
\newcommand{\inttomap}{\mathbf{M}_\mathcal{I}}

\newcommand{\degtoedge}{\Lambda}
\newcommand{\id}{\mathrm{id}}

\title{Bijective link between Chapoton's new intervals and bipartite planar maps}

\author{Wenjie Fang\thanks{Email: \href{mailto:wenjie.fang@u-pem.fr}{wenjie.fang@u-pem.fr}. Some of the work was performed during a postdoc at TU Graz financed by the Austrian Science Fund (FWF) I2309 and P27290.} \\ LIGM, Univ. Gustave Eiffel, CNRS, ESIEE Paris \\ F-77454 Marne-la-Vallée, France}

\begin{document}

\maketitle

\abstract{In 2006, Chapoton defined a class of Tamari intervals called ``new intervals'' in his enumeration of Tamari intervals, and he found that these new intervals are equi-enumerated with bipartite planar maps. We present here a direct bijection between these two classes of objects using a new object called ``degree tree''. Our bijection also gives an intuitive proof of an unpublished equi-distribution result of some statistics on new intervals given by Chapoton and Fusy.}

\section{Introduction}

On classical Catalan objects, such as Dyck paths and binary trees, we can define the famous \emph{Tamari lattice}, first proposed by Dov Tamari \cite{tamari-def}. This partial order was later found woven into the fabric of other more sophisticated objects. A notable example is diagonal coinvariant spaces \cite{bergeron-preville,hopf-dream}, which have led to several generalizations of the Tamari lattice \cite{bergeron-preville,PRV}, and also incited the interest in intervals in such Tamari-like lattices. Recently, there is a surge of interest in the enumeration \cite{ch06,bousquet-fusy-preville,chatel-pons,nonsep} and the structure \cite{bernardi-bonichon,trinity,chapoton-note} of different families of Tamari-like intervals. In particular, several bijective relations were found between various families of Tamari-like intervals and planar maps \cite{bernardi-bonichon,nonsep,sticky}. The current work is a natural extension of this line of research.

In \cite{ch06}, other than counting Tamari intervals, Chapoton also introduced a subclass of Tamari intervals called \emph{new intervals}, which are irreducible elements in a grafting construction of intervals. Definitions of these objects and related statistics are postponed to the next section. The number of new intervals in the Tamari lattice of order $n \geq 2$ was given in \cite{ch06}, which equals
\[
  \frac{3 \cdot 2^{n-2} (2n-2)!}{(n-1)!(n+1)!}.
\]
This is also the number of bipartite planar maps with $n-1$ edges. Furthermore, in a more recent unpublished result of Chapoton and Fusy (see \cite{fusy-talk} for details), a symmetry in three statistics on new intervals was observed, then also proven by identifying the generating function of new intervals recording these statistics with that of bipartite planar maps recording the number of black vertices, white vertices and faces, three statistics well-known to be equi-distributed. These results strongly hint a bijective link between the two classes of objects.

\begin{figure}
  \centering
  \includegraphics[page=1,width=\textwidth]{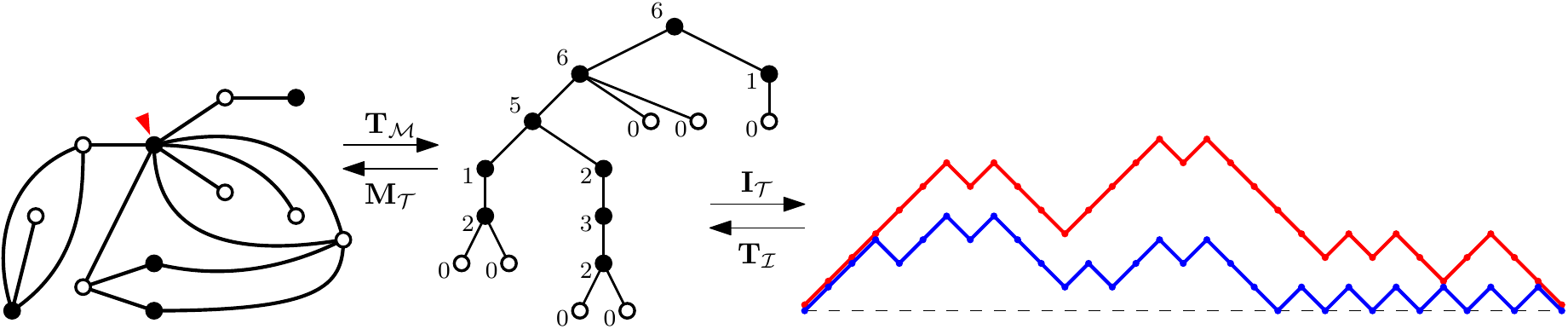}
  \caption{Our bijections between bipartite planar maps, degree trees and new intervals}
  \label{fig:main-bij}
\end{figure}

In this article, we give a direct bijection between new intervals and bipartite planar maps (see Figure~\ref{fig:main-bij}) explaining the results above. Our bijection can also be seen as a generalization of a bijection on trees given in \cite{tree-bij} in the study of random maps. We have the following theorem, with statistics defined in the next section.

\begin{thm} \label{thm:bij}
  There is a bijection $\maptoint$ from the set $\mathcal{I}_{n+1}$ of new intervals of size $n+1$ to the set $\mathcal{M}_n$ of bipartite planar maps with $n$ edges for every $n \geq 0$, with $\inttomap$ its inverse, such that, for a bipartite planar map $M$ and $I = \maptoint(M)$, which is a new interval, we have
  \begin{align*}
    \white(M) = \sdd(I), &\quad \black(M) = \sdu(I)\\
    \face(M) = 1 + \suu(I), &\quad \outdeg(M) = \rcont(I)-1.
  \end{align*}
\end{thm}

This bijection is intermediated by a new family of objects called \emph{degree trees}, and was obtained in the spirit of some previous work of the author \cite{nonsep,sticky}. Our bijection was inspired and extending another bijection given in \cite{tree-bij} between plane trees, which can be seen as bipartite planar maps.

Although the symmetry between statistics in new intervals is already known, our bijection captures this symmetry in an intuitive way, thus also opens a new door to the structural study of new intervals via bipartite maps and related objects. It is particularly interesting to see what natural involutions on bipartite maps, such as switching black and white in the coloring, induce on new intervals via our bijections.

In the rest of this article, we first define the related objects and statistics in Section~\ref{sec:prelim}. Then we show a bijection between bipartite planar maps and degree trees in Section~\ref{sec:tree-map}, then a bijection between degree trees and new intervals in Section~\ref{sec:tree-int}. We conclude by some remarks on the study of symmetries in new intervals in Section~\ref{sec:sym}.

\section{Preliminaries} \label{sec:prelim}

A \tdef{Dyck path} $P$ is a lattice path composed by up steps $u = (1,1)$ and down steps $d = (1, -1)$, starting from the origin, ending on the $x$-axis while never falling below it. A \tdef{rising contact} of $P$ is an up step of $P$ on the $x$-axis. A non-empty Dyck path has at least one rising contact, which is the first step. We can also see a Dyck path $P$ as a word in the alphabet $\{u,d\}$ such that all prefixes have more $u$ than $d$. The \tdef{size} of a Dyck path is half its length. We denote by $\mathcal{D}_n$ the set of Dyck paths of size $n$.

\begin{figure}
  \centering
  \includegraphics[page=2,width=0.7\textwidth]{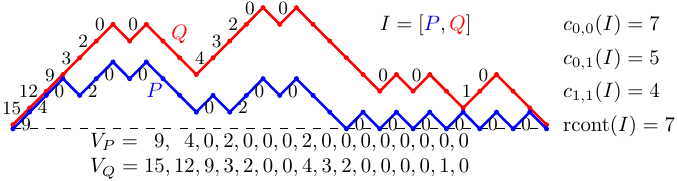}
  \caption{An example of Chapoton's new interval with bracket vectors for both paths and related statistics.}
  \label{fig:def-int}
\end{figure}

We now define the Tamari lattice, introduced in \cite{tamari-def}, as a partial order on $\mathcal{D}_n$ using a characterization in \cite{huang-tamari}. Given a Dyck path $P$ seen as a word, its $i^{\rm th}$ up step $u_i$ \tdef{matches} with a down step $d_j$ if the factor $P_i$ of $P$ strictly between $u_i$ and $d_j$ is also a Dyck path. It is clear that there is a unique match for every $u_i$. We define the \tdef{bracket vector} $V_P$ of $P$ by taking $V_P(i)$ to be the size of $P_i$. The \tdef{Tamari lattice} of order $n$ is the partial order $\preceq$ on $\mathcal{D}_n$ such that $P \preceq Q$ if and only if $V_P(i) \leq V_Q(i)$ for all $1 \leq i \leq n$. See Figure~\ref{fig:def-int} for an example. A \tdef{Tamari interval} of size $n$ can be viewed as a pair of Dyck paths $[P, Q]$ of size $n$ with $P \preceq Q$.

In \cite{ch06}, Chapoton defined a subclass of Tamari intervals called ``new intervals''. Originally defined on pairs of binary trees, this notion can also be defined on pairs of Dyck paths (see \cite{fusy-talk}). The example in Figure~\ref{fig:def-int} is also a new interval. Given a Tamari interval $[P,Q]$, it is a \tdef{new interval} if and only if the following conditions hold:
\begin{itemize}
\item[(i)] $V_Q(1) = n$;
\item[(ii)] For all $1 \leq i \leq n$, if $V_Q(i) > 0$, then $V_P(i) \leq V_Q(i+1)$.
\end{itemize}
We denote by $\mathcal{I}_n$ the set of new intervals of size $n \geq 1$.

We now define several statistics on new intervals. Given a Dyck path $P$ of size $n$, its \tdef{type} $\type(P)$ is defined as a word $w$ such that, if the $i$-th up step $u_i$ is followed by an up step in $P$, then $w_i = 1$, otherwise $w_i = 0$. Since the last up step is always followed by a down step, we have $w_n = 0$. Note that our definition here is slightly different from that in, \textit{e.g.}, \cite{nonsep}, where the last letter is not taken into account. Given a new interval $I = [P,Q] \in \mathcal{I}_n$, if $\type(P)_i=1$ and $\type(Q)_i=0$, then we have $V_P(i) > 0$ and $V_Q(i) = 0$, violating the condition for Tamari interval. Therefore, we have only three possibilities for $(\type(P)_i, \type(Q)_i)$. We define $\sdd(I)$ (resp. $\sdu(I)$ and $\suu(I)$) to be the number of indices $i$ such that $(\type(P)_i, \type(Q)_i) = (0, 0)$ (resp. $(0,1)$ and $(1,1)$). We also define $\rcont(I)$ to be the number of rising contacts of the lower path $P$ in $I = [P,Q]$. Figure~\ref{fig:def-int} also shows such statistics in the example. We define the generating function $F_{\mathcal{I}} \equiv F_{\mathcal{I}}(t,x;u,v,w)$ of new intervals as
\begin{equation}
  \label{eq:int-gf}
  F_{\mathcal{I}}(t,x;u,v,w) = \sum_{n \geq 1} t^n \sum_{I \in \mathcal{I}_n} x^{\rcont(I)-1} u^{\sdd(I)} v^{\sdu(I)} w^{\suu(I)}.
\end{equation}
We note that the power of $x$ of the contribution of a new interval $I$ is $\rcont(I)-1$.

\begin{figure}
  \centering
  \includegraphics[page=3,width=\textwidth]{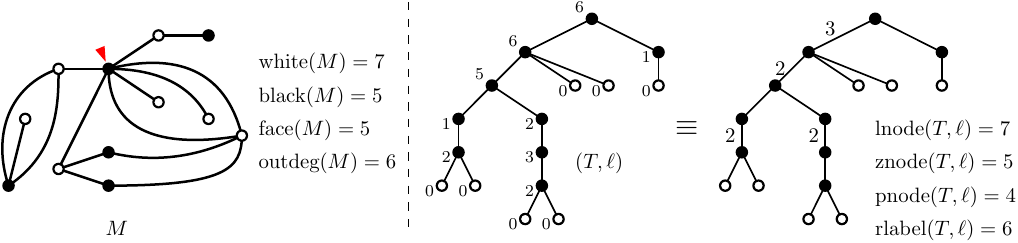}
  \caption{Left: an example of bipartite map. Right: an example of degree trees and the corresponding edge labels (zeros are omitted). Both with related statistics.}
  \label{fig:def-bip-tree}
\end{figure}

For the other side of the bijection, a \tdef{bipartite planar map} $M$ is a drawing of a bipartite graph (in which all edges link a black vertex to a white one) on the plane, defined up to continuous deformation, such that edges intersect only at their ends. Edges in $M$ cut the plane into \tdef{faces}, and the \tdef{outer face} is the infinite one. The \tdef{size} of $M$ is its number of edges. In the following, we only consider \emph{rooted} bipartite planar maps, which have a distinguished corner $c$ called the \tdef{root corner} of the outer face on a black vertex, which is called the \tdef{root vertex}. See the left part of Figure~\ref{fig:def-bip-tree} for an example. We denote by $\mathcal{M}_n$ the set of (rooted) bipartite planar maps of size $n$. We allow the bipartite planar map of size $0$, which consists of only one black vertex.

We also define some natural statistics on bipartite planar maps. For $M$ a bipartite planar map, we denote by $\black(M)$, $\white(M)$ and $\face(M)$ the number of black vertices, white vertices and faces respectively. We also denote by $\outdeg(M)$ the \emph{half-degree} of the outer face, \textit{i.e.}, half its number of corners. We take the convention that the outer face of the one-vertex map is of degree $0$.
These statistics are also illustrated in the left part of Figure~\ref{fig:def-bip-tree}. We define the generating function $F_{\mathcal{M}} \equiv F_{\mathcal{M}}(t,x;u,v,w) $ of bipartite planar maps enriched with these statistics by
\begin{equation}
  \label{eq:bip-gf}
  F_{\mathcal{M}} \equiv F_{\mathcal{M}}(t,x;u,v,w) = \sum_{n \geq 0} t^n \sum_{M \in \mathcal{M}_n} x^{\outdeg(M)} u^{\black(M)} v^{\white(M)} w^{\face(M)}.
\end{equation}
It is well known that $\black(M), \white(M), \face(M)$ are jointly equi-distributed in $\mathcal{M}_n$, meaning that $F_{\mathcal{M}}$ is symmetric in $u, v, w$. This can be seen with the bijection between bipartite maps and bicubic maps by Tutte \cite{tutte-census}, or with rotation systems of bipartite maps (see \cite[Chapter~1]{lando-zvonkin}).

To describe our bijection, we propose an intermediate class of objects called ``degree trees''. An example is given in the right part of Figure~\ref{fig:def-bip-tree}. The meaning of this name will be clear in the description of our bijection. We can also see degree trees as a variant of description trees introduced by Cori, Jacquard and Schaeffer in \cite{description-tree}. A \tdef{degree tree} is a pair $(T,\ell)$, where $T$ is a plane tree, and $\ell$ is a labeling function defined on nodes of $T$ such that
\begin{itemize}
\item If $v$ is a leaf, then $\ell(v)=0$;
\item If $v$ is an internal node with $k$ children $v_1, v_2, \ldots, v_k$, then $\ell(v) = k - a + \ell(v_1) + \ell(v_2) + \cdots + \ell(v_k)$ for some integer $a$ with $0 \leq a \leq \ell(v_1)$.
\end{itemize}
We observe that the leftmost child of a node $v$ is special when computing $\ell(v)$. This is different from the case of description trees. The size of a degree tree $(T,\ell)$ is the number of edges. We denote by $\mathcal{T}_n$ the set of degree trees $(T,\ell)$ of size $n$.

Given a degree tree $(T,\ell)$, we can replace $\ell$ by a labeling function on \emph{edges}. More precisely, for an internal node $v$, we label its leftmost descending edge by the value of $a$ used in the computation of $\ell(v)$, and all other edges by $0$. We denote this edge labeling function by $\ell_\degtoedge$. It is clear that, given $T$, the mapping $\ell \mapsto \ell_\degtoedge$ is an injection. Given $\ell_\degtoedge$, we can easily recover $\ell$ using its definition with the value $a = \ell_\degtoedge(v)$ when computing $\ell(v)$. 

We also define several natural statistics on degree trees, illustrated in Figure~\ref{fig:def-bip-tree}, using its edge labeling. Let $(T,\ell)$ be a degree tree with $\ell_\degtoedge$ the corresponding edge labeling, and $v$ a node in $T$. If $v$ is a leaf, then it is called a \tdef{leaf node}. Otherwise, let $e$ be the leftmost descending edge of $v$. If $\ell_\degtoedge(e)=0$, then $v$ is a \tdef{zero node}, otherwise it is a \tdef{positive node}. We denote by $\lnode(T,\ell)$, $\znode(T,\ell)$ and $\pnode(T,\ell)$ the number of leaf nodes, zero nodes and positive nodes in $(T,\ell)$ respectively. For $T \in \mathcal{T}_n$, we have $\lnode(T,\ell) + \znode(T,\ell)+\pnode(T,\ell) = n+1$. We also define the statistic $\rlabel$ by taking $\rlabel(T,\ell) = \ell(r)$ with $r$ the root of $T$. 

\begin{lem} \label{lem:deg-tree}
  Let $(T,\ell)$ be a degree tree, and $\ell_\degtoedge$ the corresponding edge labeling. We have
  \begin{enumerate}
  \item If $v$ has $m$ descendants, then we have $\ell(v) = m - \sum_{e \in T_v} \ell_\degtoedge(e)$, where $T_v$ is the subtree induced by $v$;
  \item $\ell(v)$ is positive, and we have $\ell(v)=0$ if and only if $v$ has no descendant.
  \end{enumerate}
\end{lem}
\begin{proof}
  The first point can be seen through induction on tree size. It holds clearly for the tree with no edge. Let $T$ be a tree of size $n$, and $v$ its root. Since the subtrees induced by each $v_i$ have sizes strictly less than $n$, by induction hypothesis, we only need to check the condition on $v$. Let $v_1, \ldots, v_k$ be the descendants of $v$, and $e_i$ the edge linking $v_i$ and $v$. From the definition of $\ell$ we have
  \[
    \ell(v) = k - \ell_\degtoedge(e_1) + \ell(v_1) + \cdots + \ell(v_k).
  \]
  To show that $\ell(v) = m -\sum_e \ell_\degtoedge(e)$, we must account for all descendants and all edges in $T_v$. However, those in one of the subtree induced by some $v_i$ are already accounted in $\ell(v_i)$. What remain are the nodes $v_1, \ldots, v_k$, which are accounted by $k$, and the edges $e_1, \ldots, e_k$, which are accounted by $-\ell_\degtoedge(e_1)$, as $\ell_\degtoedge(e_i) = 0$ for all $i > 1$. We thus conclude the induction.

  The second point can also be proved by induction on tree size. It is clearly correct when $T$ is the tree with no edge, and for the induction step, we observe that
  \[
    \ell(v) = k + (\ell(v_1) - \ell_\degtoedge(e_1)) + \ell(v_2) + \cdots + \ell(v_k) \geq k \geq 1,
  \]
  since $\ell(v_i) \geq 0$ by induction hypothesis and $0 \leq \ell_\degtoedge(e_1) \leq \ell(v_1)$ by the definition of $\ell_\degtoedge$.
\end{proof}

\section{Degree trees and bipartite maps} \label{sec:tree-map}

Our bijection from bipartite maps to new intervals is relayed by degree trees, in which the related statistics are transferred in an intuitive way. We now start by the bijection from maps to trees.

\subsection{From bipartite maps to degree trees}

It is well known that plane trees with $n$ nodes in which $k$ of them are leaves are counted by Narayana numbers (\textit{cf.} \cite{drmota-tree}). In \cite{tree-bij}, Janson and Stef\'ansson described a bijection between such plane trees and plane trees with $n$ nodes in which $k$ of them are of even depth, providing yet another interpretation of Narayana numbers. We now introduce a bijection between bipartite planar maps and degree trees, which can be seen as a generalization of the bijection in \cite{tree-bij}.

\begin{figure}
  \centering
  \includegraphics[page=4,width=0.9\textwidth]{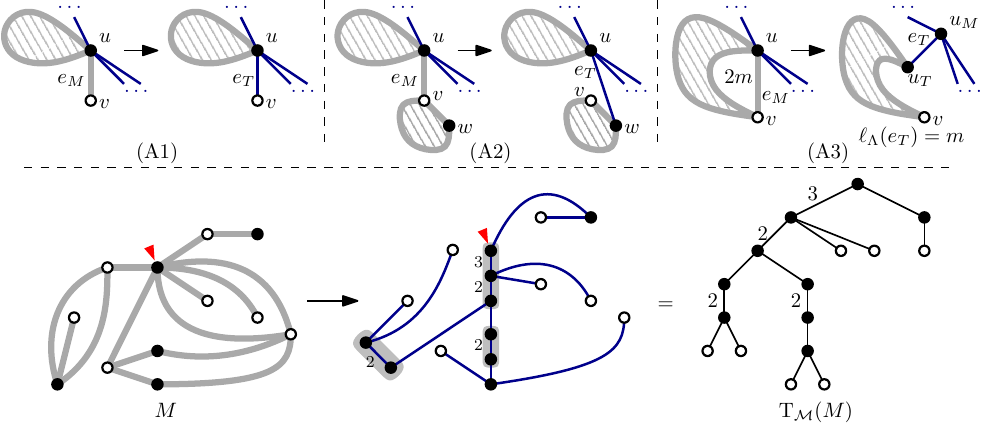}
  \caption{Cases in the advance step of $\maptotree$ and an example of the bijection $\maptotree$. Nodes in the same shaded pack come from the same vertex in the map.}
  \label{fig:map-bij}
\end{figure}

We first define a transformation $\maptotree$ from $\mathcal{M}_n$ to $\mathcal{T}_n$ for all $n$. Let $M \in \mathcal{M}_n$. If $n=0$, we define $\maptotree(M)$ to be the tree with one node. Otherwise, we perform the following exploration procedure to obtain a tree $T$ with a labeling $\ell_\degtoedge$ on its edges. In this procedure, we distinguish edges in $M$, which will be deleted one by one, and edges in $T$ that we add. We start from the root vertex, with the edge next to the root corner in clockwise order as the pending edge. Suppose that the current vertex is $u$ and the pending edge is $e_M$, which is always in $M$. We repeat two steps, \tdef{advance} and \tdef{prepare}, until termination. Roughly, in the advance step we modify edges in $M$ and $T$ and update the current vertex and the pending edge, and then in the prepare step we fix potential problems. The advance step comes in the following cases illustrated in Figure~\ref{fig:map-bij}:
\begin{itemize}
\item[\textbf{(A1)}] If $e_M$ is a bridge to a vertex $v$ of degree $1$, then we delete $e_M$ in $M$ and add $e_T = e_M$ in $T$. The new current vertex is $u' = u$, and we define $\ell_\degtoedge(e_T) = 0$.
\item[\textbf{(A2)}] If $e_M$ is a bridge to a vertex $v$ of degree at least $2$, let $e_1$ be the edge adjacent to $v$ next to $e_M$ in clockwise order, and $w$ the other end of $e_1$. We draw a new edge $e_T$ in $T$ from $u$ to $w$ such that $e_M, e_1, e_T$ form a face with $u,v,w$ in counter-clockwise order. The next current vertex is $u' = w$. We delete $e_M$, and define $\ell_\degtoedge(e_T)=0$.
\item[\textbf{(A3)}] If $e_M$ is not a bridge, we split $u$ into $u_M$ and $u_T$, with $u_T$ taking all edges in $T$ and $u_M$ taking the rest. We add a new edge $e_T$ in $T$ from $u_M$ to $u_T$. Since $e_M$ is not a bridge, by planarity, it is between the outer face and a face of degree $2m$ with $m>0$. We define $\ell_\degtoedge(e_T) = m$ and delete $e_M$. The next current vertex is $u' = u_T$.
\end{itemize}
In the prepare step, let $u'$ be the new current vertex, which is adjacent to the new edge $e_T$. The next pending edge is the next remaining edge in $M$ starting from $e_T$ in the clockwise order around $u'$. If no such edge exists, we backtrack in the tree $T$ until finding a vertex $u''$ with such an edge $e_M''$, and we set $u''$ as the current vertex, and $e_M''$ the pending edge. If no such vertex exists, the procedure terminates, and we shall obtain a tree $T$ with an edge label function $\ell_\degtoedge$. We define $\maptotree(M)$ as the degree tree $(T, \ell)$, with $\ell$ the node labeling corresponding to $\ell_\degtoedge$. See Figure~\ref{fig:map-bij} for an example of $\maptotree$. The bijection in \cite{tree-bij} is simply $\maptotree$ applied to a plane tree, where Case~(A3) never applies, and the degree tree $(T,\ell)$ obtained has $\ell_\degtoedge=0$ for all edges.

We now prove that $\maptotree(M)$ is well-defined. We start by describing the structure of the map in intermediate steps. The \tdef{leftmost branch} of a tree is the path starting from the root node and taking the leftmost descending edge at each node till a leaf.

\begin{lem} \label{lem:map2tree-struct}
  Let $M \in \mathcal{M}_n$ and $T = \maptotree(M)$. Let $M^+_i$ be the map after the $i$-th prepare step, with $u_i$ the current vertex and $e_i$ the pending edge. We denote by $T_i$ the partially constructed $T$ in $M^+_i$ , and $M_i$ that of the remaining of $M$. Clearly $T_i$ and $M_i$ form a partition of edges in $M^+_i$.

  For every $i$, $T_i$ is a tree, and $M^+_i$ is $T_i$ with connected components of $M_i$ attached to the left of nodes on the leftmost branch of $T_i$, one component to only one vertex, with $u_i$ the deepest such vertex and $e_i$ its first edge in $M_i$ in clockwise order from the leftmost branch of $T_i$.
\end{lem}
\begin{proof}
  We proceed by induction on $i$. The case $i=0$ is trivial. We now suppose that the induction hypothesis holds for $i$, and we prove that it also holds for $i+1$. Suppose that the component of $M_i$ attached to $u_i$ is $M_{i,*}$, then $e_i$ is in $M_{i,*}$. For the $(i+1)$-st advance step, we have three possibilities.
  \begin{itemize}
  \item \textbf{Case (A1)}: $e_i$ links $u_i$ to a node $v_i$ of degree $1$. The advance step then turns $e_i$ into an edge in $T_{i+1}$. It is clear that $T_{i+1}$ is also a tree, and other components of $M_i$ are still in $M_{i+1}$ and attached to the same vertices, except $M_{i,*}$, which becomes empty if only $e_i$ is in it, or is turned into $M_{i+1,*}$ with $e_i$ deleted otherwise. In the latter case, since $v_i$ was of degree $1$, the deletion of $e_i$ does not disconnect $M_{i+1,*}$, thus $M_{i+1,*}$ is still attached to $u_i$. Either way, all components of $M_{i+1}$ are still attached to $T_{i+1}$ on the leftmost branch. Then in the prepare step, either $M_{i+1,*}$ is not empty, and we have $u_{i+1} = u_i$, with $e_{i+1}$ the next edge in clockwise order of $e_i$, or it is empty, and we backtrack on the leftmost branch until finding a vertex with a component of $M_{i+1}$ attached, which is also the last one in the preorder of $T_{i+1}$, and $e_{i+1}$ is the next edge in the clockwise order of the last backtracking edge. Therefore, by induction hypothesis, $e_{i+1}$ is also the first one in $M_{i+1}$ starting from any edge of $u_{i+1}$ in $T_{i+1}$. 
  \item \textbf{Case (A2)}: $e_i$ links $u_i$ to a node $v_i$ of degree at least $2$, and $e_i$ is a bridge in $M_i$, thus also in $M_{i,*}$. The removal of $e_i$ breaks $M_{i,*}$ into two parts, $M_{i+1,1}$ attached to $u_i$, and $M_{i+1,2}$ containing $v_i$. Let $e_{T,i}$ be the edge added to $T_{i+1}$ in the advance step, linking $u_i$ to a node $w_i$. By construction, $w_i$ is in $M_{i+1,2}$, therefore not in $T_i$ by induction hypothesis. Thus, $T_{i+1}$ is a tree, and the newly separated component $M_{i+1,2}$ is attached to $T_{i+1}$ by $w_i$. All other components of $M_i$ remains in $M_{i+1}$ and attached to $T_{i+1}$. Then in the prepare step, since $M_{i+1,2}$ is not empty, we have $u_{i+1} = w_i$, and $e_{i+1}$ the first edge of $w_i$ in $M_{i+1,2}$ in clockwise order, starting from $e_i$ linking $w_i$ to its parent $u_i$. 
  \item \textbf{Case (A3)}: $e_i$ is not a bridge in $M_i$. The remaining $M_{i+1,*}$ of $M_{i,*}$ after the removal of $e_i$ is still connected. Let $e_{T,i}$ be the edge added to $T_{i+1}$ in the advance step, linking $u_i$ to a node $w_i$. By construction, $M_{i+1,*}$ is attached to $w_i$. We verify the conditions on $u_{i+1}$ and $e_{i+1}$ with the same reasoning as in Case~(A2). 
  \end{itemize}
  As the induction hypothesis is valid in all cases, we conclude the proof.
\end{proof}

We now prove that trees obtained in $\maptotree$ are degree trees.

\begin{prop} \label{prop:map2tree}
  Given $M \in \mathcal{M}_n$ a bipartite map of size $n$, the tree $(T,\ell) = \maptotree(M)$ is a degree tree of size $n$.
\end{prop}
\begin{proof}
  From Lemma~\ref{lem:map2tree-struct}, we know that the whole procedure of $\maptotree$ does not stop before consuming all $n$ edges in $M$, and $T$ is a tree. Therefore, $T$ is a tree of size $n$.

  Let $\ell_\degtoedge$ be the edge labeling obtained in the procedure of $\maptotree$. The labels in $\ell_\degtoedge$ are all positive by construction. We also observe that $\ell_\degtoedge(e) > 0$ for an edge $e \in T$ implies that $e$ links a node $u$ to its leftmost child, as only Case~(A3) has the possibility of $\ell_\degtoedge(e) > 0$, and the new edge $e_T$ added in that case becomes the leftmost descending edge of $u$ after the duplication. We now only need to prove that the node labeling $\ell$ corresponding to $\ell_\degtoedge$ satisfies the conditions of degree trees.

  We now define a labeling $\ell'$ on nodes of $T$. By Lemma~\ref{lem:map2tree-struct}, the first time a node $u$ is explored on $T$, there is a component of some remaining edges in $M$ attached to $u$, which is itself a planar map. We denote by $M_u$ this planar map. We define $\ell'(u)$ to be half of the degree of the outer face of $M_u$. We now prove that $\ell(u) = \ell'(u)$ by induction on the size of the subtree induced by $u$. For the base case, $u$ is a leaf, and $\ell(u) = 0 = \ell'(u)$. When $u$ is an internal node with children $u_1, \ldots, u_k$ from left to right, by induction hypothesis, we have $\ell(u_i) = \ell'(u_i)$ for all $i$. Now, for $i \geq 2$, the node $u_i$ is produced by Case~(A1)~or~(A2), thus are linked by bridges to $u$ in $M_u$. The contribution of such $u_i$ to $\ell'(u)$ is thus $\ell'(u_i)+1$. For $u_1$, by checking all cases, its contribution to $\ell'(u)$ is $\ell'(u_1)+1-\ell_\degtoedge(e_1)$, where $e_1$ is the edge between $u$ and $u_1$. The only case that needs attention is Case~(A3), where a face of degree $2\ell_\degtoedge(e_1)$ is merged with the outer face by the removal of $e_1$, increasing the degree of the outer face by $2\ell_\degtoedge(e_1)-2$. Therefore, the degree of the outer face of the part attached to $u$ leading to $u_1$ before the exploration of $u_1$ is the correct value $2\ell'(u_1)+2-2\ell_\degtoedge(e_1)$. We thus have
  \[
    \ell'(u) = \sum_{i=1}^k (\ell'(u_1)+1) - \ell_\degtoedge(e_1) = k - \ell_\degtoedge(e_1) + \sum_{i=1}^k \ell(u_1) = \ell(u).
  \]
  We thus conclude by induction that $\ell = \ell'$. Then, since the degree of the outer face of a planar bipartite map is at least $2$, we have $\ell(u_1)-\ell_\degtoedge(e_1) \geq 0$ for each edge $e_1$ from a node $u$ to its first child $u_1$. Hence, $(T,\ell)$ satisfies the conditions of degree trees.
\end{proof}

The transformation $\maptotree$ transfers some statistics from $\mathcal{M}_n$ to $\mathcal{T}_n$ as follows.

\begin{prop} \label{prop:map-bij-stat}
  Given $M \in \mathcal{M}_n$, let $(T,\ell) = \maptotree(M)$. We have
  \begin{align*}
    \white(M) = \lnode(T,\ell), &\quad \black(M) = \znode(T,\ell), \\
    \face(M) = 1 + \pnode(T,\ell), &\quad \outdeg(M) = \rlabel(T, \ell).
  \end{align*}
\end{prop}
\begin{proof}
  Since in $\maptotree$ we only walk on black vertices, all leaves in $T$ are from white vertices, which are never split. Hence $\white(M) = \lnode(T,\ell)$. Then at each occurrence of Case~(A3), we lost a face but gain a positive node in $T$, thus $\face(M) = 1 + \pnode(T,\ell)$, with $1$ for the outer face. Now for $\black(M) = \znode(T,\ell)$, we note that a new black vertex in $M$ is reached only in Case~(A2), which leads to a zero edge. For $\outdeg(M)$, we notice $\outdeg(M) = n - \sum_{f} \deg(f)/2$, summing over all internal faces $f$ of $M$. However, by the bijection, we have $\sum_{f} \deg(f)/2 = \sum_{e \in T} \ell_\degtoedge(e)$, and we conclude by Lemma~\ref{lem:deg-tree}(1) applied to the root.
\end{proof}

\subsection{From degree trees to bipartite maps}

\begin{figure}
  \centering
  \includegraphics[page=5,width=0.9\textwidth]{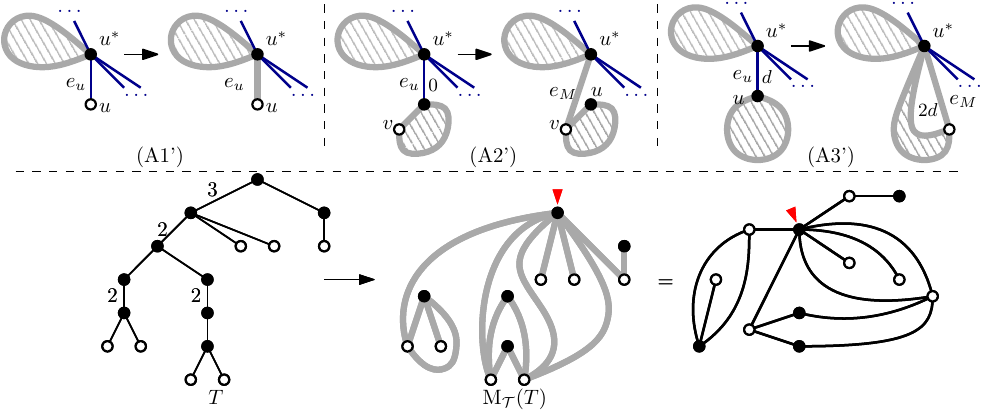}
  \caption{Cases in the procedure of $\treetomap$, and an example of $\treetomap$}
  \label{fig:map-bij-inv}
\end{figure}

We now define a transformation $\treetomap$ from $\mathcal{T}_n$ to $\mathcal{M}_n$, which is precisely the inverse of $\maptotree$. Let $(T,\ell) \in \mathcal{T}_n$ and $\ell_\degtoedge = \degtoedge(\ell)$. We now perform the following procedure that deals with nodes in $T$ in postorder (\textit{i.e.}, first visit the subtrees induced by children from left to right, then the parent). For each node $u$, let $u^*$ be its parent and $e_u$ the edge between $u$ and $u^*$. By construction, when we deal with $u$, its induced subtree has already been dealt with, transformed into a bipartite planar map $M_u$ attached to $u$. We have three cases, illustrated in Figure~\ref{fig:map-bij-inv}.
\begin{itemize}
\item \textbf{Case~(A1')}: If $u$ is a leaf, then we delete $e_u$ from $T$ and add it to $M$.
\item \textbf{Case~(A2')}: If $u$ is not a leaf but $\ell_\degtoedge(e_u) = 0$, let $e'$ be the edge next to $e_u$ around $u$ in counterclockwise order, and $v$ the other end of $e'$. As $M_u$ is bipartite, $v \neq u$. We add a new edge $e_M$ from $u^*$ to $v$ such that the triangle formed by $e_u, e', e_M$ has vertices $u^*, u, v$ in clockwise order, without any edge inside. We then delete $e_u$.
\item \textbf{Case~(A3')}: If $\ell_\degtoedge(e_u) > 0$, let $d$ be the degree of the outer face of $M_u$. If $2 \ell_\degtoedge(e_u) > d$, then the procedure fails. Otherwise, we start from the corner of $M_u$ to the right of $e_u$ and walk clockwise along edges for $2 \ell_\degtoedge(e_u) - 1$ times to another corner, and we connect the two corners by a new edge $e_M$ in $M$, making a new face of degree $2\ell_\degtoedge(e_u)$. The component remains planar and bipartite. We finish by contracting $e_u$.
\end{itemize}
In the end, we obtain a planar bipartite map $M$ with the same root corner as $T$. We define $\treetomap(T,\ell) = M$. We see that (A1'), (A2') and (A3') are exactly the opposite of (A1), (A2), (A3) in the definition of $\maptotree$.

We first show that the procedure above never fails, thus $\treetomap$ is always well-defined. It follows easily that we always have bipartite planar maps from $\treetomap$.

\begin{prop} \label{prop:tree2map-defined}
  Given $(T,\ell)$ a degree tree, for a node $u \in T$, let $M_u$ be the map obtained in the procedure of $\treetomap(T,\ell)$ from the subtree $T_u$ induced by $u$. Then the degree of the outer face of $M_u$ is $2\ell(u)$, and the procedure never fails.
\end{prop}
\begin{proof}
  We use induction on the size of the subtree $T_u$. It clearly holds when $u$ is a leaf. Suppose that $u$ is an internal node. Let $u_1, \ldots, u_k$ be its children from left to right. Since every edge $e_i$ linking $u_i$ to $u$ must be in Case~(A1')~or~(A2') for $i \geq 2$, the contribution of the part $M_{u_i}$ to the degree of the outer face is $2 + 2\ell(u_i)$ by induction hypothesis. If $e_1$ linking $u_1$ to $u$ is also a bridge, then the contribution is $2 + 2\ell(u_i)$. Otherwise, we are in Case~(A3'), in which we create a new face of degree $2\ell_\degtoedge(e)$, where $\ell_\degtoedge$ is the corresponding edge labeling. We never fail in this case, since by the definition of $\degtoedge$, we have $0 \leq \ell_\degtoedge(e) \leq \ell(u_1)$. Therefore, $M_{u_1}$ has an outer face of degree $2\ell(u_1) + 2 - 2\ell_\degtoedge(e)$. The degree of the outer face of $M_u$ is thus
  \[
    2\ell(u_1) + 2 - 2\ell_\degtoedge(e) + \sum_{i=2}^k (2 + 2\ell(u_i)) = 2\ell(u).
  \]
  We thus conclude the induction.
\end{proof}

\begin{prop} \label{prop:tree2map-valid}
  For $(T,\ell)$ a degree tree, $M = \treetomap(T,\ell)$ is a bipartite planar map.
\end{prop}
\begin{proof}
  Planarity is easily checked through the definition of $\treetomap$. Faces in $M$ are only created in Case~(A3'), thus all of even degree. Since $M$ is planar, every cycle of edges can be seen as a gluing of faces, which are all of even degree. Therefore, the cycle obtained is always of even length, meaning that $M$ is bipartite.
\end{proof}

It is also clear that $\treetomap$ is the inverse of $\maptotree$.

\begin{prop} \label{prop:map-tree-bij}
  The transformation $\maptotree$ is a bijection from $\mathcal{M}_n$ to $\mathcal{T}_n$, with $\treetomap$ its inverse.  
\end{prop}
\begin{proof}
  By Proposition~\ref{prop:map2tree}, we only need to prove that $\maptotree \circ \treetomap = \id_{\mathcal{T}}$ and $\treetomap \circ \maptotree = \id_{\mathcal{M}}$.

  For $\maptotree \circ \treetomap = \id_{\mathcal{T}}$, it is clear that the operations in cases of $\treetomap$ are reverted by those in $\maptotree$, and by Lemma~\ref{lem:map2tree-struct}, the degree tree is constructed node by node in reverse postorder in $\maptotree$. We thus have $\maptotree \circ \treetomap = \id_{\mathcal{T}}$.

  To show that $\treetomap \circ \maptotree = \id_{\mathcal{M}}$, we only need to check that they are applied exactly in the reverse order, and there is only one possibility for reversing operations in each case of $\maptotree$. The first point is again ensured by Lemma~\ref{lem:map2tree-struct}. For the second point, the only case to check is Case~(A3). To revert operation in this case, we need to create a new face of given degree by cutting the outer face with an edge. By planarity, there is only one way to proceed, which is that of Case~(A3) in $\treetomap$. We thus conclude that $\treetomap$ is indeed the inverse of $\maptotree$, and they are all bijections.
\end{proof}

\section{Degree trees and new intervals} \label{sec:tree-int}

We now present the bijective link between degree trees and new intervals, which also gives a combinatorial explanation of the conditions of new intervals in terms of trees.

\subsection{From degree trees to new intervals}

\begin{figure}
  \centering
  \includegraphics[page=6,width=\textwidth]{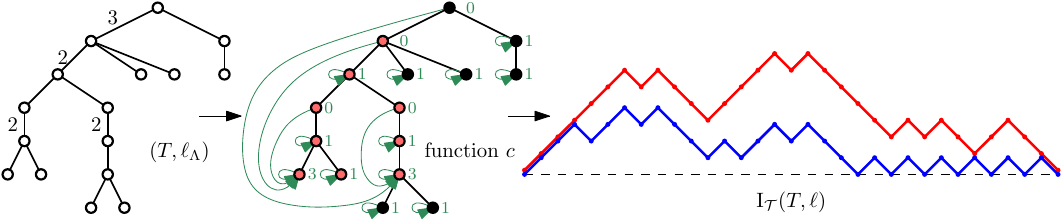}
  \caption{Example of the bijection $\treetoint$ on a degree tree represented by its edge labeling. The middle shows the certificate of each node.}
  \label{fig:int-bij}
\end{figure}

Given $(T,\ell) \in \mathcal{T}_n$, let $\ell_\degtoedge$ be the corresponding edge labeling. We define a transformation $\treetoint$ by constructing a pair of Dyck paths $[P,Q]$ from $(T,\ell)$. We take $Q = uQ'd$, where $Q'$ comes from the classical bijection between plane trees and Dyck paths by doing a traversal of $T$ in \emph{preorder} (parent first, then subtrees from left to right), recording the evolution of depth. For $P$, we first assign to every node a \tdef{certificate}, and we define a \tdef{certificate function} $c$ on $T$ as in \cite{sticky,nonsep}. We process all nodes in $T$ in the \emph{reverse} preorder, initially colored black. At the step for a node $v$, if $v$ is a leaf, then its certificate is itself. Otherwise, let $e$ be the leftmost descending edge of $v$. We then visit nodes after $v$ in preorder, and color each visited black node by red. We stop at the node $w$ just before the $(\ell_\degtoedge(e)+1)$-st black node, and the certificate of $v$ is $w$. When $\ell_\degtoedge(e)=0$, we take $w=v$. Now, we take $c(w)$ to be the number of nodes with $w$ as certificate. With the function $c$, the path $P$ is given by concatenation of $ud^{c(v)}$ for all nodes $v$ in preorder. We then define $\treetoint(T,\ell) = [P,Q]$. An example of $\treetoint$ is given in Figure~\ref{fig:int-bij}.

To prove that $\treetoint(T,\ell)$ is a new interval, we start by some properties of certificates.

\begin{lem} \label{lem:tree-cert}
  Let $(T,\ell)$ be a degree tree of size $n$ and $\ell_\degtoedge$ the corresponding edge labeling. For a node $v \in T$, let $w$ be the certificate of $v$. Then either $w = v$, or $w$ is a descendant of $v$ in the leftmost subtree $T_*$ of $v$. In the latter case, $w$ is not the last node of $T_*$ in preorder.
\end{lem}
\begin{proof}
  Let $v_1, v_2, \ldots, v_{n+1}$ be the nodes in $T$ in preorder. We prove our statement for all $v_i$ by reverse induction on $i$. It is clear that the last node $v_n$ in preorder is a leaf, hence its certificate is itself. The base case is thus valid.

  For the induction step, suppose that all $v_j$'s with $i < j \leq n$ satisfy the induction hypothesis. If $v_i$ is a leaf, then the induction hypothesis holds for $i$. We now suppose that $v_i$ has at least one child. Let $T_*$ the subtree induced by the left-most child $v_*$ of $v_i$, and $e_*$ the edge linking $v_*$ to $v_i$. If $v_*$ is a leaf, then $\ell_\degtoedge(e_*)=0$ and the induction hypothesis is clearly correct. We suppose that $v_*$ is not a leaf. We consider the coloring just before the step for $v_i$. Since nodes in $T_*$ come after $v_i$ in the preorder, their processing only changes color of nodes in $T_*$ by induction hypothesis. Therefore, there are $\sum_{e \in T_*} \ell_\degtoedge(e)$ red nodes in $T_*$. By Lemma~\ref{lem:deg-tree}(1), there are thus $(\ell(v_*) + 1)$ black nodes in $T_*$, where the extra $1$ accounts for $v_*$ itself, which is never red after its process step. Since $\ell_\degtoedge(e) \leq \ell(v_*) + 1$, the $(\ell(v_*) + 1)$-st black node starting from $v_*$ must be in $T_*$. Hence, the certificate of $v_i$ is either $v_i$ or in $T_*$, and cannot be the last node in $T_*$. We thus conclude the induction.
\end{proof}

\begin{lem} \label{lem:cert-parenthesis}
  Let $(T,\ell)$ be a degree tree, and $v, v'$ two distinct nodes in $T$ with $w, w'$ their certificates respectively. Suppose that $v$ precedes $v'$ in the preorder. Then $w$ cannot be strictly between $v'$ and $w'$ in the preorder. Furthermore, if $v' \neq w'$, then $w \neq v'$.
\end{lem}
\begin{proof}
  We only need to consider the case $v \neq w$ and $v' \neq w'$, as other cases are trivial. In the coloring process, since $v$ precedes $v'$ in the preorder, $v'$ is treated before $v$. By construction, in the coloring process, after the step for $v'$, the nodes between $v'$ to $w'$ (excluding $v'$ but including $w'$) are all colored red. Therefore, in the process step for $v$, the visit will not stop strictly between $v'$ and $w'$, nor at $v'$, as such a stop requires a succeeding black node. Hence, $w$ is not strictly between $v'$ and $w'$, and $w \neq v'$.
\end{proof}

Note that in the lemma above, we can have $w = v'$ when $v' = w'$.

\begin{prop} \label{prop:tree2int}
  Let $(T,\ell) \in \mathcal{T}_n$ be a degree tree of size $n$. The pair of Dyck paths $[P,Q] = \treetoint(T,\ell)$ is a new interval in $\mathcal{I}_{n+1}$.
\end{prop}
\begin{proof}
  Let $v_1, v_2, \ldots, v_{n+1}$ be the nodes in $T$ (including the root) in preorder, and $T_i$ the subtree induced by $v_i$ for $1 \leq i \leq n+1$. We now prove that both $P$ and $Q$ are Dyck paths, with a combinatorial interpretation of their bracket vector $V_P$ and $V_Q$. From the construction of $Q$, it is clear that $Q$ is a Dyck path, and we have $V_Q(i) = |T_i|$, where $|T_i|$ is the size of $T_i$ (\textit{i.e.}, the number of edges).

  For $P$, from the construction of $P$ and Lemma~\ref{lem:tree-cert}, a node that gives an up step never comes after its certificate that gives a down step, meaning that there are at least as many up steps as down steps in any prefix of $P$, making it a Dyck path. To compute $V_P(i)$, we consider $v_i$, its certificate $w_i$, and the subword $P'$ of $P$ that comes from the nodes from $v_i$ to $w_i$ (both $v_i$ and $w_i$ included). If $v_i$ is a leaf or $v_i = w_i$, it is clear that $P' = ud^{c(v_i)}$ and $V_P(i)=0$. Otherwise, we consider a node $v_j$ strictly between $v_i$ and $w_i$ in the preorder of $T$, in which case we can write $P' = u d^{c(v_i)} P'' u d^{c(w_i)}$. Firstly, let $w_j$ be the certificate of $v_j$, then by Lemma~\ref{lem:cert-parenthesis}, $w_j$ cannot come strictly after $w_i$. Thus in $P'$ there are more down steps than up steps. Secondly, by Lemma~\ref{lem:cert-parenthesis}, no node has $v_i$ as certificate, implying that $c(v_i)=0$. Thirdly, also by Lemma~\ref{lem:cert-parenthesis}, if $v_j$ is a certificate of a node, then this node must be strictly between $v_i$ and $v_j$, already contributing an up step to $P''$. Therefore, in any prefix of $P''$, there are at least the same number of up steps than down steps. We then have the $i$-th up step in $P$ generated by $v_i$ matches with one of the down steps in $P'$ (by the first point), but not those in $P''$ or induced by $v_i$ itself (by the second and the third point), therefore it matches with a down step generated by $w_i$. Since $v_{i+1}$ is the first child of $v_i$. By Lemma~\ref{lem:tree-cert}, $w_i$ is in the subtree induced by $v_{i+1}$, but not the last node, implying $V_P(i) \leq |T_{i+1}|$.

  We now compare $V_P$ and $V_Q$. It is clear that $V_Q(1) = n$. If $V_Q(i) = 0$, then $v_i$ is a leaf, and we have $V_P(i) = 0 \leq V_Q(i)$. If $V_Q(i) > 0$, then $v_i$ has descendants, and we have $V_P(i) \leq |T_{i+1}| = V_Q(i+1) < V_Q(i)$ in this case. Therefore, the pair $[P,Q]$ is not only a Tamari interval, but also a new interval. It is clear from the construction of $P$ and $Q$ that they are Dyck paths of size $n+1$.
\end{proof}

We also have the following property of the new interval obtained from a given degree tree via $\treetoint$.

\begin{prop} \label{prop:tree-edge-contact}
  For a degree tree $(T,\ell)$ with $\ell_\degtoedge$ the corresponding edge labeling, let $I = [P,Q] = \treetoint(T,\ell)$. For an internal node $v \in T$, let $e$ be the edge linking $v$ to its leftmost child $v'$, and $r = \ell_\degtoedge(e)$. Let $P_v$ be the subpath of $P$ strictly between the up step contributed by $v$ in $P$ and its matching down step. Then the number of rising contacts in $P_v$ as a Dyck path is $r$.
\end{prop}
\begin{proof}
  Let $w$ be the certificate of $v$. The subpath $P_v$ comes from the contributions of nodes from $v'$ to $w$, while deleting extra down steps from $w$ due to potentially other nodes preceding $v$ in preorder taking up $w$ as certificate.

  By Lemma~\ref{lem:cert-parenthesis}, no node preceding $v$ in preorder has its certificate strictly between $v$ and $w$, and the certificate of nodes from $v'$ to $w$ cannot be strictly after $w$ in the preorder. Therefore, $P_v$ is totally determined by the relation of certificates for nodes from $v'$ to $w$, which is known when the coloring process gets $v$ treated. In that step, exactly $r$ black nodes are colored red, denoted by $v_1, v_2, \ldots, v_r$ in the preorder. Let $w_1, \ldots, w_r$ be their certificates respectively.

  First we prove that, for $1 \leq i \leq r$, the subpath of $P_v$ contributed by nodes from $v_i$ to $w_i$, denoted by $P_i$, is a Dyck path with one rising contacts. This is again due to Lemma~\ref{lem:cert-parenthesis}, making the certificates of nodes strictly between $v_i$ and $w_i$ to be between $v_i$ and $w_i$ (can be equal to $w_i$). Thus $P_i$ has the same number of up steps and down steps. Since the up step from a node always comes before the down step from its certificate, $P_i$ is a Dyck path. There is no other rising contact of $P_i$, because the up step from $v_i$ is matched by the last down step from $w_i$.

  Now, clearly we have $v_1=v'$, as $v'$ is the node next to $v$ in preorder, thus treated in the coloring process just before $v$, but the treatment always leave $v_1$ black. Now, at the step of $v_1$ in the coloring process, $w_1$ is the red node just before a black node in preorder. This black node cannot come after $v_2$, as it would entail $v_2$ being red in the step for $v$, but not before $v_2$ either, as it would still be black in the step for $v$, violating the definition of $v_2$. The same argument applies to all $v_i$, thus the next node of $w_i$ in preorder is $v_{i+1}$ for $1 \leq i \leq r-1$. We now consider $w_r$. The node $v_+$ next to $w_r$ in preorder must be black at the step for $v_r$, and remains black through all treatments for nodes till $v_1$. Therefore, $v_+$ must come strictly after $w$, and we can only have $w = w_r$. We thus conclude that every node from $v'$ to $w$ is between some pair of $v_i$ and $w_i$. Therefore, we can write $P_v = P_1 \cdots P_r$, and we conclude that the number of rising contacts in $P_v$ is indeed $r$.
\end{proof}

\subsection{From new intervals to degree trees}

\begin{figure}
  \centering
  \includegraphics[page=7,width=\textwidth]{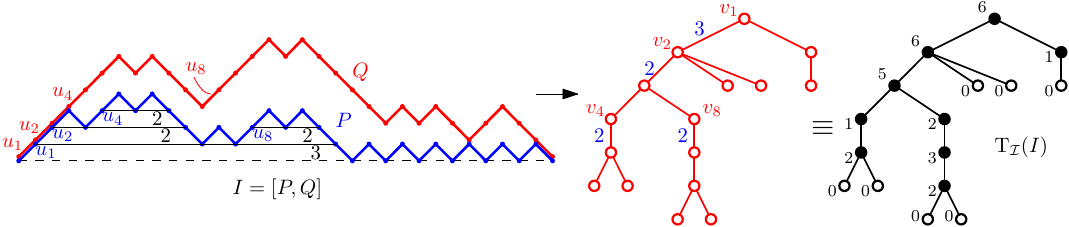}
  \caption{Example of the bijection $\inttotree$ on a new interval $I = [P,Q]$}
  \label{fig:int-bij-inv}
\end{figure}

We now define a transformation $\inttotree$ for the reverse direction. Let $I = [P,Q] \in \mathcal{I}_{n+1}$ be a new interval. Since $V_Q(1) = n$, we can write $Q = uQ'd$. We first construct a plane tree $T$ of size $n$ from $Q'$ using again the classical bijection. Now, let $v_1, \ldots, v_{n+1}$ be the nodes of $T$ in preorder. We note that $V_Q(i)$ is the size of the subtree induced by $v_i$, which is equal to the number of descendants of $v_i$. We now define the edge labeling $\ell_\degtoedge$ of $T$. If $e$ is the left-most descending edge of $v_i$, then we take $\ell_\degtoedge(e)$ the number of rising contacts in $P_i$, where $P_i$ is the subpath of $P$ strictly between the $i$-th up step and its matching down step. Otherwise, we take $\ell_\degtoedge(e) = 0$. We define $\inttotree(I) = (T, \ell)$, with $\ell$ the node labeling corresponding to $\ell_\degtoedge$. An example of $\inttotree$ is given in Figure~\ref{fig:int-bij-inv}. We first show that $(T,\ell)$ is indeed a degree tree.

\begin{prop} \label{prop:int2tree}
  Let $I = [P,Q] \in \mathcal{I}_{n+1}$, then $(T,\ell) = \inttotree(I)$ is a degree tree of size $n$.
\end{prop}
\begin{proof}
  Let $\ell_\degtoedge$ be the edge labeling obtained when applying $\inttotree$ to $I$. We start by the following property of $\ell_\degtoedge$. Suppose that $e'$ is an edge linking the $j$-th node $v_j$ in $T$ to its leftmost child $v_{j+1}$, and $T_{j+1}$ is the subtree induced by $v_{j+1}$. We know that $\ell_\degtoedge(e)$ is the number of rising contacts in $P_j$, where $P_j$ is the subpath of $P$ strictly between the $j$-the up step and its matching down step. In other words, $\ell_\degtoedge(e')$ is the number of up steps in $P_j$ that starts at the same height ($y$-coordinate) as the upper end of the $j$-th up step in $P$. Since in this case we have $V_Q(j) > 0$ as $v_j$ is not a leaf, by the condition of new intervals, we have $V_P(j) \leq V_Q(j+1)$. Since up steps in $P_j$ comes from descendants of $v_j$, and $V_Q(j+1)$ is the number of descendants of $v_{j+1}$, which are the first descendants of $v_j$ in preorder, we conclude that all up steps in $P$ contributing to $\ell_\degtoedge(e')$ are from nodes in $T_{j+1}$, but not the last one in preorder.

  From the construction, it is clear that the sizes match, and we only need to show that, for any edge $e$ linking an internal node $v$ to its leftmost child $v_*$, we have $\ell_\degtoedge(e) \leq \ell(v_*)$. Let $m$ be the number of descendants of $v_*$, and $T_*$ is the subtree induced by $v_*$. The property above means that nodes whose up steps contributed to $\ell_\degtoedge(e)$ or $\ell_\degtoedge(e')$ for any $e' \in T_*$ must be in $T_*$, but not the last one in preorder. It is clear that every up step can only contribute to $\ell_\degtoedge(e')$ for at most one $e'$. We thus have
  \[
    m \geq \ell_\degtoedge(e) + \sum_{e' \in T_*} \ell_\degtoedge(e').
  \]
  We deduce $\ell_\degtoedge(e) \leq \ell(v_*)$ using the same argument as for Lemma~\ref{lem:deg-tree}(1).
\end{proof}

Some natural statistics are transferred from new intervals to degree trees via $\inttotree$.

\begin{prop} \label{prop:int-bij-stat}
  Given $I = [P,Q] \in \mathcal{I}_{n+1}$, let $(T,\ell) = \inttotree(I)$. We have
  \begin{align*}
    \sdd(I) = \lnode(T,\ell), &\quad \sdu(I) = \znode(T,\ell), \\
    \suu(I) = \pnode(T,\ell), &\quad \rcont(I) = 1 + \rlabel(T,\ell).
  \end{align*}
\end{prop}
\begin{proof}
  Let $v_i$ be the $i$-th node of $T$ in preorder. By the definition of $\inttotree$, the node $v_i$ is a leaf if and only if $\type(Q)_i=0$. Hence, $\sdd(I) = \lnode(T,\ell)$. Moreover, if $v_i$ is an internal node, then $\type(P)_i=0$ if and only if $\ell_\degtoedge(e_i) = 0$, where $e_i$ is the leftmost descending edge of $v_i$, and $\ell_\degtoedge$ the edge labeling corresponding to $\ell$. We thus conclude for $\sdu(I) = \znode(T,\ell)$ and $\suu(I) = \pnode(T,\ell)$. For $\rcont(I)$, we observe that rise contacts come from up steps not contributing to the edge labeling $\ell_\degtoedge$, meaning that $\rcont(I) = n + 1 - \sum_{e \in T} \ell_\degtoedge(e)$. By applying Lemma~\ref{lem:deg-tree}(1) to the root, we have $\rlabel(T,\ell) = n - \sum_{e \in T} \ell_\degtoedge(e)$, therefore $\rcont(I) = 1 + \rlabel(T,\ell)$.
\end{proof}

Using Proposition~\ref{prop:tree-edge-contact}, we check that $\treetoint$ and $\inttotree$ are bijections.

\begin{prop} \label{lem:int-tree-bij}
  For any $n \geq 0$, the transformation $\treetoint$ is a bijection from $\mathcal{T}_n$ to $\mathcal{I}_{n+1}$, with $\inttotree$ its inverse.
\end{prop}
\begin{proof}
  By Propositions~\ref{prop:tree2int}~and~\ref{prop:int2tree}, we only need $\inttotree \circ \treetoint = \id_{\mathcal{T}}$ and $\treetoint \circ \inttotree = \id_{\mathcal{I}}$.

  For $\inttotree \circ \treetoint = \id_{\mathcal{T}}$, let $(T,\ell) \in \mathcal{T}_n$ and $I = [P,Q] = \treetoint(T,\ell)$. Now we consider $(T', \ell') = \inttotree(I)$. It is clear from the definition of $\inttotree$ and $\treetoint$ that $T = T'$. We now show that $\ell = \ell'$, which is equivalent to $\ell_\degtoedge = \ell'_\degtoedge$, where $\ell_\degtoedge$ (resp. $\ell'_\degtoedge$ is the edge labeling corresponding to $\ell$ (resp. $\ell'$). Let $e$ be an edge in $T$. We only need to consider the case where $e$ links a node $v$ to its leftmost child $v'$. Suppose that $v$ is the $i$-th node in the preorder of $T$. Let $P_i$ be the subpath of $P$ between the $i$-th up step and its matching down step. Now by Proposition~\ref{prop:tree-edge-contact} and the definition of $\inttotree$, the number of rising contacts in $P_i$ is equal to both $\ell_\degtoedge(e)$ and $\ell'_\degtoedge(e)$, making $\ell_\degtoedge(e) = \ell'_\degtoedge(e)$, thus $\ell = \ell'$. We conclude that $\inttotree \circ \treetoint = \id_{\mathcal{T}}$.
  
  For $\treetoint \circ \inttotree = \id_{\mathcal{I}}$, let $I = [P,Q] \in \mathcal{I}_{n+1}$ and $(T,\ell) = \inttotree(I)$. We take $\ell_\degtoedge$ the edge labeling corresponding to $\ell$. Now we consider $I' = [P',Q'] = \treetoint(I)$. Again, it is clear that $Q=Q'$, and we only need to show that $P = P'$. For $1 \leq i \leq n+1$, let $P_i$ (resp. $P'_i$) be the subpath of $P$ (resp. $P'$) strictly between the $i$-th up step and its matching down step, and $e_i$ the edge linking the $i$-th node in the preorder of $T$ to its leftmost child. By the definition of $\inttotree$ and Proposition~\ref{prop:tree-edge-contact}, there are $\ell_\degtoedge(e_i)$ rising contacts in both $P_i$ and $P'_i$ for every $i$. However, suppose that $P$ (resp. $P'$) leads to a plane tree $T_P$ (resp. $T_{P'}$) via the classical bijection. Since the number of rising contacts in $P_i$ (resp. $P'_i$) is the degree of the $(i+1)$-st node in the preorder of $T_P$ (resp. $T_{P'}$), we know that the degrees of nodes in $T_P$ and $T_{P'}$ in preorder are the same. This leads to $T_P = T_{P'}$, meaning that $P = P'$. We thus conclude that $\treetoint \circ \inttotree = \id_{\mathcal{I}}$.
\end{proof}

\section{Symmetries and structure} \label{sec:sym}

With the bijections in Section~\ref{sec:tree-map}~and~\ref{sec:tree-int}, we construct the following bijections between new intervals and bipartite maps, which is our main result.

\begin{proof}[Proof of Theorem~\ref{thm:bij}]
  We take $\maptoint = \treetoint \circ \maptotree$ and $\inttomap = \treetomap \circ \inttotree$. Their validity is from Proposition~\ref{prop:map-tree-bij}~and~\ref{lem:int-tree-bij}. The equalities of statistics come from Proposition~\ref{prop:map-bij-stat}~and~\ref{prop:int-bij-stat}.
\end{proof}

The symmetry between the statistics $\white$, $\black$ and $\face$ on bipartite maps is then transferred to new intervals.

\begin{coro} \label{coro:stats-sym}
  The generating functions $F_\mathcal{I}$ and $F_\mathcal{M}$ are related by
  \[
    t F_\mathcal{M} = w F_\mathcal{I}.
  \]
  In particular, the series $w F_\mathcal{I}$ is symmetric in $u, v, w$.
\end{coro}
\begin{proof}
  The equality is a direct translation of Theorem~\ref{thm:bij} in generating functions. The symmetry of $w F_\mathcal{I}$ comes from that of $F_\mathcal{M}$.
\end{proof}

\begin{rmk}
  For $M \in \mathcal{M}$, let $D$ be the multiset of half-degrees of internal faces of $M$. From the definition of $\maptotree$, the multiset $D$ is also the multiset of non-zero edge labels of $\maptotree(D)$. Now, let $[P,Q] = \maptoint(M)$. From the definition of $\inttotree$, the multiset of non-zero labels of $\maptotree(D)$ is also that of the number of rising contacts of subpaths of $P$ between matching steps. We can thus refine Corollaory~\ref{coro:stats-sym} by this multiset $D$. Such refinement is particularly interesting in the domain of maps. We can enrich $F_\mathcal{M}$ by an infinity of variables $(p_k)_{k \geq 1}$, with $p_k$ marking internal faces of half-degree $k$. Such enriched version is particularly nice and has deep link with factorization of the symmetric group and other objects. See \cite{cf-bipartite} for more details. It would be interesting to see how results on refined enumeration of bipartite maps can be transferred to new intervals.
\end{rmk}

As mentioned before, the symmetry of $\sdd, \sdu, \suu$ in new intervals was already known to Chapoton and Fusy, and a proof relying on generating functions was outlined in \cite{fusy-talk}, which makes use of recursive decompositions of new intervals \cite[Lemma~7.1]{ch06} and bipartite planar maps. Our bijective proof can be seen as direct version of this recursive proof, in the sense that $\inttotree$ and $\maptotree$ are canonical bijections of these recursive decompositions. Details will be given in a follow-up article.

Since there are bijections for bipartite maps that permute black vertices, white vertices and faces arbitrarily, there should be an isomorphic symmetry structure hidden in new intervals via bijections. If we regard new intervals as pairs of binary trees, it is easy to see that there is an involution consisting of exchanging the two trees in the pair while taking their mirror images. This involution exchanges the statistics $\sdd$ and $\suu$, corresponding to $\white$ and $\face$ in bipartite planar maps. The structural study of these symmetries under our bijections is the subject of a follow-up article.

Furthermore, these is another class of combinatorial objects called $\beta$-(0,1) trees, which are description trees for bicubic planar maps in bijection with bipartite maps \cite{description-tree,desc-tree-tutte}. An involution on these trees is given in \cite{involution-beta01}, which may be related to symmetries we mentioned above. 

However, as a precaution for all structural study, we should note that our bijections are subjected to various choices taken in their definitions. For instance, in the definition of the bijection $\maptotree$ from bipartite maps to degree trees, in the case~(A2), if we fix an integer $k \geq 0$, and construct the new edge $e_T$ in the partial tree by connecting to the $k$-th black corner on the outer face in clockwise order, and changing the definition of $\treetomap$ accordingly, we will have a bijection parameterized by $k$, which is different for all $k$. There is thus an infinity of bijections compatible with all results in this article. Therefore, it is possible that the bijections defined here may not preserve some wanted structure between related objects, but a similar bijection does.

As pointed out by an anonymous reviewer, the notion of degree tree bears similarities to that of ``grafting trees'' defined in \cite{rise-contact}, which is in turn closely related to description trees of type $(1, 1)$ in \cite{desc-tree-tutte} and closed flows on forests \cite{chapoton-chatel-pons, sticky}. Since our degree tree can be seen as a special case of description trees of type $(1,1)$, it would be interesting to see how the bijections extend to the general Tamari intervals and corresponding maps. For instance, the anonymous reviewer also observed that, elements of our bijection from degree trees to new intervals can be used to construct a direct bijection from grafting trees to general Tamari intervals.

\section*{Acknowledgment}

The author thanks \'Eric Fusy for inspiring discussions, especially about the recursive decomposition of related objects from his own discussion with Frédéric Chapoton. The author also thanks Philippe Biane and Samuele Giraudo for proofreading the conference version of this article. Moreover, the author also thanks the anonymous reviewers' precious and interesting comments.

\bibliographystyle{alpha}
\bibliography{bip-tam}

\end{document}